\documentclass{amsart}
\usepackage{amssymb, amscd}
\newcommand{\ie}{{\em i.e., }}

\newcommand{\CIc}{{\mathcal C}^{\infty}_{\text{c}}}
\newcommand{\pa}{\partial}

\newcommand{\maD}{\mathcal D}

\newcommand{\maQ}{\mathcal Q}

\newcommand{\CC}{\mathbb C}
\newcommand{\NN}{\mathbb N}
\newcommand{\RR}{\mathbb R}

\newcommand{\ZZ}{\mathbb Z}

\newcommand{\pairing}[2]{\langle #1, #2 \rangle}
\newcommand{\innerp}[2]{( #1, #2 )}

\newcommand\normH[2]{\|#1\|_{H^{#2}(\Omega)}}

\newcommand{\loc}{\mathrm{loc}}

\newcommand{\supp}{\operatorname{supp}}
\newcommand{\diam}{\operatorname{diam}}
\newcommand\ha{\frac{1}{2}}

\newcommand{\CORR}{\ }

\newtheorem{theorem}{Theorem}[section]
\newtheorem{proposition}[theorem]{Proposition}
\newtheorem{corollary}[theorem]{Corollary}

\newtheorem{lemma}[theorem]{Lemma}

\theoremstyle{definition}
\newtheorem{definition}[theorem]{Definition}
\theoremstyle{remark}
\newtheorem{remark}[theorem]{Remark}
\newtheorem{example}[theorem]{Example}

\author[I Babu\v{s}ka]{Ivo Babu\v{s}ka} \address{The University of
Texas at Austin, Institute for Computational Engineering and
Sciences, Austin, TX 78712--0027}
\email{Babuska@ticam.utexas.edu}

\author[V. Nistor]{Victor Nistor} \address{Pennsylvania State
University, Math. Dept., University Park, PA 16802}
\email{nistor@math.psu.edu}

\thanks{I. Babu\v{s}ka was partially supported by NSF Grant DMS
0341982. V. Nistor was partially supported by NSF Grants DMS 991981
and 0200808.  Manuscripts available from {\bf
http:{\scriptsize//}www.math.psu.edu{\scriptsize/}nistor{\scriptsize/}}.}

\begin{document}

\dedicatory\date

\title[Distributions]{Interior numerical approximation
of boundary value problems with a distributional data}

\begin{abstract}
We study the approximation properties of a harmonic function $u
\in H\sp{1-k}(\Omega)$, $k > 0$, on relatively compact sub-domain
$A$ of $\Omega$, using the Generalized Finite Element Method. For
smooth, bounded domains $\Omega$, we obtain that the
GFEM--approximation $u_S$ satisfies $\|u - u_S\|_{H\sp{1}(A)} \le
C h^{\gamma}\|u\|_{H\sp{1-k}(\Omega)}$, where $h$ is the typical
size of the ``elements'' defining the GFEM--space $S$ and $\gamma
\ge 0 $ is such that the local approximation spaces contain all
polynomials of degree $k + \gamma + 1$. The main technical result
is an extension of the classical super-approximation results of
Nitsche and Schatz \cite{NitscheSchatz72} and, especially,
\cite{NitscheSchatz74}. It turns out that, in addition to the
usual ``energy'' Sobolev spaces $H^1$, one must use also the
negative order Sobolev spaces $H\sp{-l}$, $l \ge 0$, which are
defined by duality and contain the distributional boundary data.
\end{abstract}

\maketitle \tableofcontents

\section*{Introduction}

To motivate our results, let us consider the boundary value
problem
\begin{equation}
\begin{cases}
\label{eq.BVP}
    \,\Delta u = 0 & \text{ on } \;\Omega,\\
    \,\frac{\pa u}{\pa \nu} = g & \text{ on }
    \Gamma := \partial \Omega,
\end{cases}
\end{equation}
where $\Omega$ is a {\em smooth, bounded} domain in $\RR^n$, $\pa
\Omega$ is the boundary of $\Omega$ and $g \in H^{r-3/2}(\pa
\Omega)$, $r \in \RR$. The case $r > 3/2$ was extensively studied.
Here we are interested mainly in the case $r \le 3/2$. We are
looking for a solution $u \in H^{r}(\Omega)$.

For $r > 3/2$, the boundary values (or {\em traces}) $u\vert_{\pa
\Omega}$ and $\pa_{\nu}u\vert_{\pa \Omega}$ are defined
classically, because the restriction to the boundary extends by
continuity to a map $H^{r}(\Omega) \ni u \to \pa_\nu u \in
H^{r-3/2}(\pa \Omega)$, see \cite{Evans, Taylor1} for example. For
$r \le 3/2$, this is no longer true, but then one takes advantage
of the fact that $u$ satisfies an elliptic equation, so it is
still possible to define $\pa_\nu u \in H^{r-3/2}(\pa \Omega)$,
see \cite{hor63, LionsMagenes1}. We can assume, without loss of
generality, that $\Omega$ is connected. It is not difficult to
prove that a solution $u$ of Equation \eqref{eq.BVP} exists for
any $g$ such that $\pairing{g}{1} = 0$ and that this solution
satisfies
\begin{equation}\label{eq.gl.est}
    \|u\|_{H^{r}(\Omega)} \le C
    \|g\|_{H\sp{r-3/2}(\pa \Omega)},
\end{equation}
with a constant $C$ that may depend on $r$ but is independent of
$g$. This result will be discussed in detail in
\cite{BabuNistor2}, where more references will be given. (This
result will be used in this paper only as a motivation for our
work.)

Recall that $A \Subset \Omega$ means that $A \subset \Omega$ is
bounded and $\pa A$ and $\pa \Omega$ are disjoint. If $A \Subset
\Omega$ is an open subset, then $u$ satisfying \eqref{eq.BVP} will
be smooth on $A$, regardless of what $r$ is, and
\begin{equation}\label{eq.int.est}
    \|u\|_{H^{m}(A)} \le C \|u\|_{H\sp{r}(\pa \Omega)},
\end{equation}
with a constant $C$ that depends on $A$, $\Omega$, $r$, and $m$,
but is independent of $u$ satisfying $\Delta u = 0$. An important
problem, with potential practical applications, is to approximate
on $A$ the solution $u$ of Equation \eqref{eq.BVP}.

In this paper, we prove several results on the approximation of
the solution $u$ on subsets $A \Subset \Omega$. The main result of
the paper, Theorem \ref{thm.main}, then gives in particular that
the solution $u$ of the boundary value problem \eqref{eq.BVP}
satisfies
\begin{equation}\label{eq.thm.main}
    \|u - u_S\|_{H\sp{1}(A)} \le C h\sp{\gamma}
    \normH{u}{1-k},
\end{equation}
where $u_S \in S$ is the Galerkin approximation and $S$ is the
Generalized Finite Element Space associated to a partition of
unity $\{\phi_j\}$ subordinated to a covering $\{\omega_j\}$ of
$\Omega$ satisfying Assumptions A--D of Section \ref{sec.GFEM},
provided that our local approximation spaces contain all of
polynomials of degree $1 + \gamma + k$, $\gamma \ge 0$.
We stress that our results require not just the energy Sobolev
space $H\sp{1}$, but also negative order Sobolev spaces
$H\sp{-l}$, defined in this paper as the duals of $H\sp{l}$, $l
\in \ZZ_+$. One of the main reasons for the need to considere the
negative order Sobolev spaces is that the solution $u$ is not in
$H\sp{1}(\Omega)$ but in $H\sp{1-k}(\Omega)$, in general.
Moreover, even if we approximate the boundary data $g$ and the
solution $u$ with functions in $H\sp{1}$, then the norm on a
negative order Sobolev space will still have to appear in the
estimate of the error.

Here is now a brief description of the contents of the paper. {\em
We continue to assume that $\Omega$ is bounded and connected, but
we no longer assume that $\Omega$ is smooth, except when
explicitly mentioned.} In Section \ref{sec.prel} we set up the
notation and we establish our conventions on Sobolev spaces.
Section \ref{sec.GFEM} contains a quick review of the necessary
definitions involving the Generalized Finite Element Method (GFEM)
and their variants used in this paper \cite{BabuBaOs,
BabuskaMelenk1, BabuskaMelenk2}. This section also contains the
assumptions that we make on our covering and partition of unity
used to define the GFEM--space $S$. The space $S$ will contain our
approximate solution to the boundary value problem. The following
section, Section \ref{sec.IE}, contains the calculations necessary
to establish our interior estimates. \CORR Our approach follows,
to a certain extend, that in the article of Nitsche and Schatz
\cite{NitscheSchatz74}, relying also from Wahlbin's survey
article \cite{Wahlbin91}. The main differences between our paper
and \cite{NitscheSchatz74, Wahlbin91} are due mostly to the fact
that several assumptions from these papers are not fully satisfied
in our approach. As in these articles, the main step is a
super-approximation property, Proposition \ref{prop.lemma.3.1}.
The proof in \cite{NitscheSchatz74, Wahlbin91} cannot be used to
obtain Proposition \ref{prop.lemma.3.1} because the property
``$\pa^{\alpha} u = 0$ if $|\alpha|$ is large,'' is not satisfied
in general for $u \in S$. For the results of Section \ref{sec.AS},
we assume that $\Omega$ is smooth. Then we extend the definitions
of the Galerkin approximation $u_S \in S$ and of the form $B(w,
v):= \int_{\Omega} \nabla u(x) \cdot \nabla v(x) dx$ to the case
when $v \in H^{k-1}(\Omega)$ is arbitrary and $u \in
H^{1-k}(\Omega)$ can be written as $u = u_1 + u_2$, where $\Delta
u_1 = 0$ in distributions sense and $u_2 \in H^1(\Omega)$. In what
follows, $1-k$ will play the role of $r$ above. Several estimates
for $u$ and its approximation $u_S$ are established in this
section, including the main theorem, Theorem \ref{thm.main} (whose
main conclusion was repeated in Equation \eqref{eq.thm.main}
above). The last section, Section \ref{sec.Ver} contains a proof
that, for a domain with piecewise $C^1$-boundary, we can construct
a family of partitions of unity with typical size of supports $h
\to 0$ that satisfies Assumptions A--D for a fixed choice of the
other values of the parameters (\ie of $A$, $B$, $C_j$, $\kappa$,
$m$, and $\sigma$). For this construction we assume that the local
approximation spaces are $\Psi_j = \maQ_k$, the space of
polynomials of degree at most $k$. In particular, the various
assumptions made in the results proved in the preceding sections,
are satisfied for this family of partition of unity, and hence our
results are not empty for domains with piecewise smooth boundary.
By contrast, it is not possible to find a family of partitions of
unity as above for domains with cusps, see Remark
\ref{remark.Ver}. For suitable $g$, we plan to perform some
concrete numerical simulations in a future paper.

We shall write $x := y$ if $x$ is {\em defined} by $y$.

The second named author thanks G. Grubb, A. Schatz, and L. Wahlbin
for some useful references. \CORR A. Schatz has also made some	
useful comments on an earlier version of the manuscript, for which
we are greatful.

\section{Preliminaries\label{sec.prel}}

We begin by fixing the notation and terminology.

\subsection{Preliminary notation} We denote by $\RR$ the set of
real numbers and by $\CC := \{a + b\imath, a, b \in \RR\}$ the set
of complex numbers. Also, $\NN = \{1, 2, \ldots \}$ and $\ZZ_+ =
\{0\} \cup \NN$. Let $x \cdot y := x_1 y_1 + x_2 y_2 + \ldots +
x_n y_n$ be the inner product of two vectors $x, y \in \RR^n$. We
shall denote by
\begin{equation}
    \hat{f}(y) := \int_{\RR^n} e^{-iy \cdot x} f(x) dx
\end{equation}
the Fourier transform of $f$, as usual. By $L^2(\Omega)$ we shall
denote the space of square integrable functions $f : \Omega \to
\CC$, for any domain $\Omega \subset \RR^n$.

Let $s \ge 0$. Then $H^s(\RR^n)$ is the space of functions $f \in
L^2(\RR^n)$ such that
\begin{equation}\label{eq.norm.s}
    \|f\|_{H^s(\RR^n)}^2 := (2\pi)^{-n}\int_{\RR^n}
    (1 + |y|^2)^{s} |\hat{f}(y)|^2 d y < \infty.
\end{equation}
Here $|y| = \sqrt{y_1^2 + \ldots + y_n^2}$, if $y = (y_1, \ldots,
y_n) \in \RR^n$. Let $\Omega \subset \RR^n$, then $H^s(\Omega)$,
$s \ge 0$, denotes the restrictions to $\Omega$ of functions $f
\in H^s(\RR^n)$, that is,
\begin{equation}
    H^s(\Omega) := \{f\vert_{\Omega},\ f \in H^s(\RR^n)\}.
\end{equation}
The norm on $H^s(\Omega)$ is then $\normH{h}{s} := \inf
\|f\|_{H^s(\RR^n)}$, the infimum being taken over all functions $u
\in H^s(\RR^n)$ such that $f \vert_{\Omega} = h$. When $s$ is a
positive integer and $\Omega$ is a nice domain (Lipschitz, for
example), the norm $\normH{v}{s}$ is equivalent to the usual norm
\begin{equation}
    \|h\|_s^2 = \sum_{|\alpha| \le s}
    \|\pa^\alpha h\|^2_{L^2(\Omega)},
\end{equation}
where $\alpha = (\alpha_1, \ldots, \alpha_n) \in \ZZ_+^n$,
$|\alpha| := \alpha_1 + \alpha_2 + \ldots + \alpha_n$, and
$\pa^\alpha := \pa_1^{\alpha_1} \pa_2^{\alpha_2}\ldots
\pa_n^{\alpha_n}$, as usual. The space $H^{s}_0(\Omega)$ is
defined as the closure of $\CIc(\Omega)$ in $H^{s}_0(\Omega)$.

\subsection{Distributions}
Let $B_R(0)$ denote the open ball of radius $R$ centered at the
origin. Also, let $\CIc(\RR^n)$ be the set of infinitely
differentiable, complex valued functions that vanish outside the
ball $B_R(0)$, for some large $R > 0$.  A linear map $u :
\CIc(\RR^n) \to \CC$ is called a {\em distribution on $\RR^n$}
\cite{GelfandShilov, hor1, Taylor1} if, for any $R > 0$, there
exists $m \in \ZZ_+= \{0, 1, 2, \ldots\}$ and $C > 0$ such that
\begin{equation}\label{eq.ineg.distr}
    |u(\phi)| \le C \sum_{|\alpha| \le m} \|\pa^\alpha
    \phi \|_{L^\infty}, \quad \text{if}\;\; \phi \in \CIc(\RR^n)
    \text{ and }
    \phi = 0 \;\; \text{outside} \;\; B_R(0).
\end{equation}
This definition does not exclude the case when larger and larger
values of $m$ have to be chosen as $R \to \infty$, and in fact
this situation actually occurs in specific examples. The set of
distributions on $\RR^n$ will be denoted $\maD'(\RR^n)$.

We now fix more notation and terminology. If $f$ is a function,
then the closure of the set $\{f \neq 0\}$ is called the {\em
support of $f$} and will be denoted $\supp(f)$. We shall also
write $\pairing{u}{\phi} := u(\phi)$ for the value of the
distribution $u$ on the function $\phi \in \CIc(\RR^n)$.
Therefore, any $\phi \in \CIc(\Omega)$ has compact support. The
support of a distribution $u$ is the smallest closed set $F$ such
that $\pairing{u}{\phi} = 0$ for any $\phi \in \CIc(\RR^n
\smallsetminus F)$.
\smallskip

Here are some examples of distributions.

\begin{example}\ If $f $ is a measurable function on $\RR^n$ that is
integrable on any closed ball in $\RR^n$ (\ie it is {\em locally
integrable}, or $f \in L^1_{\loc}(\RR^n)$), then we can define
\begin{equation}
    \pairing{f}{\phi} := \int_{\RR^n} f(x) \phi(x) dx,
\end{equation}
for any $\phi \in \CIc(\RR^n)$. Thus any $f \in L^1_{\loc}(\RR^n)$
defines a distribution on $\RR^n$, that is, $L^1_{\loc}(\RR^n)
\subset \maD'(\RR^n)$. In this situation we shall say that our
distribution $f$ is, in fact, a locally integrable function, or,
that our distribution is defined by a locally integrable function.
\end{example}

The following two examples are relevant for the discussion of
concentrated loads and moments.

\begin{example}\ The Dirac measure at $a\in \RR^n$ is the distribution
$\delta_a$ defined by
\begin{equation*}
    \pairing{\delta_a}{\phi} := \phi(a).
\end{equation*}
An explicit calculation shows that $\delta_a \in H^{-n/2 -
\epsilon}(\RR^n)$ and $\|\delta_a\|_{H^{-n/2 - \epsilon}(\RR^n)}
\to \infty$ as $\epsilon \to 0$ as $\epsilon^{1/2}$.
\end{example}

\begin{example}\ The derivatives $\pa^\alpha u$ of a distribution
$u$ are defined by
\begin{equation*}
    \pairing{\pa^\alpha u}{\phi} := (-1)^{|\alpha|}
    \pairing{u}{\pa^\alpha \phi}.
\end{equation*}
\end{example}

We now define the negative index Sobolev spaces. Thus, the space
$H^{-s}(\RR^n)$, $s \ge 0$, consists of all the distributions $u
\in \maD'(\RR^n)$ such that there exists $h \in L^1_{\loc}(\RR^n)$
satisfying
\begin{equation*}
    \pairing{u}{\hat{f}} = \int_{\RR^n} h(y) f(y) dy, \quad
    \text{for all }\; f \in \CIc(\RR^n)
\end{equation*}
and
\begin{equation}\label{eq.norm.s2}
    \|u\|_{-s}^2 := (2\pi)^{-n}\int_{\RR^n}
    (1 + |y|^2)^{s} |h(y)|^2 d y < \infty.
\end{equation}
The reader will recognize that the condition above is analogous to
the condition of Equation \eqref{eq.norm.s}. The difference is
that now we allow $s$ to take on negative values as well.

The following alternative definition of the negative order Sobolev
spaces will be useful later on. One first checks directly using
the Fourier inversion formula together with Plancherel's formula
that
\begin{equation}\label{eq.norm.s3}
    \|u\|_{-s}^2 := \inf_{\phi} \frac{1}{\|\phi\|_{s}} \left |
    \int_{\RR^n} u(x) \phi(x) dx \right |, \quad \text{where }\;
    u, \phi \in \CIc(\RR^n),\ v \neq 0.
\end{equation}
Then $H^{-s}(\RR^n)$ can be defined as the completion of
$\CIc(\RR^n)$ in the norm given by the formula of Equation
\eqref{eq.norm.s3}. In particular, we obtain that $H^{-s}(\RR^n)$
is canonically isomorphic to the dual of $H^s(\RR^n)$.

Similarly, for any open subset $\Omega \subset \RR^n$, we define
the space $H^{-s}(\Omega)$, $s \ge 0$, as the dual of
$H^{s}(\Omega)$.
\begin{equation*}
    \|u\|_{H\sp{-s}(\Omega)} = \sup \frac{|(u,
    \phi)|}{\|\phi\|_{H\sp{s}(\Omega)}}, \quad
    0 \neq \phi \in H\sp{s}(\Omega).
\end{equation*}

Our definition of negative order Sobolev spaces by duality follows
\cite{Ciarlet91, RoitbergBook, Schechter60}, for example. Note
however, that the negative order Sobolev spaces are often also
defined by restriction from $\RR^n$, as in \cite{Evans,
LionsMagenes1, Taylor1}, for example. The space of restrictions to
$\Omega$ of distributions in $H^{-s}(\RR^n)$ is the dual of
$H^{s}_0(\Omega)$, and will be denoted $H_{0}^{-s}(\Omega)$. The
spaces $H_{0}^{-s}(\Omega) = H_0^{s}(\Omega)^*$, $s \ge 0$, will
also be used below.

When $\Omega = \RR^n$, these two approaches yield the same spaces,
but for general $\Omega$ they may lead to different ``negative
order'' Sobolev spaces.

\section{The Generalized Finite Element Method\label{sec.GFEM}}
We now recall a few basic facts about the Generalized Finite
Element Method \cite{BabuBaOs, BabuskaMelenk1, BabuskaMelenk2}.
This method is especially convenient since it provides finite
element spaces with high regularity. Most of the results of this
section work for a general bounded open set $\Omega$, except the
application to the boundary value problems, Subsection
\ref{ssec.DS}.

\subsection{Basic facts}
We shall need the following slight generalization of a definition
from \cite{BabuskaMelenk1, BabuskaMelenk2}:

\begin{definition}\label{def.MCG}\
Let $\Omega \subset \RR^n$ be an open set and $\{\omega_j\}$ be an
open cover of $\Omega$ with no $\kappa + 1$ of the sets $\omega_j$
having a non-empty intersection. Also, let $\{\phi_j\}$ be a
partition of unity consisting of Lipschitz functions and
subordinated to the covering $\{\omega_j\}$ (\ie $\supp \phi_j
\subset \omega_j$). If
\begin{equation}\label{eq.2.10}
    \|\pa^\alpha \phi_j\|_{L^\infty(\Omega)}
    \le C_k/(\diam \omega_j)\sp{k},
    \quad k = |\alpha| \le m,
\end{equation}
then $\{\omega_j\}$ is called a {\em $(\kappa, C_0, C_1, \ldots, C_m)$
partition of unity}.
\end{definition}

Assume also that linear subspaces $\Psi_j \in H^m(\omega_j)$ are
given. These subspaces will be called {\em local approximation
spaces} and are used to define the space
\begin{equation}\label{eq.2.14}
    S = S_{GFEM} := \Big\{ \sum_j \phi_j v_j,\;  v_j \in \Psi_j
    \Big\},
\end{equation}
which will be called the {\em GFEM--space}.

A basic approximation property of the GFEM--spaces is the
following Theorem from \cite{BabuskaMelenk1}.

\begin{theorem}[Babu\v{s}ka-Melenk]\label{thm.2.1}\
We use the notations and definitions of Definition \ref{def.MCG}
and after. Let $\{\phi_j\}$ be a $(\kappa, C_0, C_1)$ partition of
unity. Also, let $v_j \in \Psi_j$, $u_{ap} := \sum_{j} \phi_j v_j
\in S$, and $d_j = \diam \omega_j$, the diameter of $\omega_j$.
Then
\begin{equation}\label{eq.2.16}
\begin{gathered}
    \|u - u_{ap} \|_{L^2(\Omega)}^2 \le \kappa C_{0}^2
    \sum_{j} \|u - v_j\|_{L^2(\omega_j)}^2 \quad \text{and}
    \\ \|\nabla (u - u_{ap})\|_{L^2(\Omega)}^2 \le 2\kappa
    \sum_{j} \Big(\frac{\,C_1^2 \|u -
    v_j\|_{L^2(\omega_j)}^2}{(d_j)^2}
    \,+\, C_0^2 \|\nabla (u - v_j)\|_{L^2(\omega_j)}^2 \Big)
\end{gathered}
\end{equation}
\end{theorem}

For our method, we shall need to make some additional assumptions
on the local approximation spaces $\Psi_j$, on the covering
$\{\omega_j\}$, and on the partition of unity $\{\phi_j\}$. We
shall denote by $d_j$ the diameter of $\omega_j$. Recall that we
have assumed that $\Omega$ is connected, the general case being
immediately reduced to this one. Unless otherwise mentioned, we
shall make the following assumptions, for some fixed values of the
parameters $A$, $B$, $C_j$, $\kappa$, $m$, $\sigma$, and $h$.
\smallskip

It follows that
\begin{equation}\label{eq.res.ABCD}
    \sum_{\omega_j^* \subset \Omega} \phi_j = 1
\quad \text{on } \Omega.
\end{equation}
In turn, Equation \eqref{eq.res.ABCD} implies the second condition
in Assumption A. The covering $\{\omega_j\}$ will satisfy the
following geometric assumption.
\smallskip

\noindent {\bf Assumption A.}\ The sets $\omega_j$ are convex of
diameters $ d_j \le h \le 1$ and there exists a ball $\omega^*_j$
of diameter greater or equal $\sigma h$ whose closure is contained
in $\omega_j$, for all $j$. Moreover, $\omega^*_j \subset \Omega$
if $\phi_j$ is not identically zero on $\Omega$.
\smallskip

This is a non-trivial assumption, see Section \ref{sec.Ver},
Remark \ref{remark.Ver}. The partition of unity will satisfy the
following condition.
\smallskip

\noindent {\bf Assumption B.}\ The family $\{\phi_j\}$ is a
$(\kappa, C_0, C_1, \ldots, C_m)$ partition of unity such that
\begin{equation}\label{eq.2.7a}
    \phi_j(x) = 1 \quad \text{if } x \in \omega^*_j.
\end{equation}

It follows from our assumptions that the sets $\omega^*_j$ must be
disjoint.  Also, we have that
\begin{equation}\label{eq.new.inverse}
    \|\phi_j\|_{H\sp{l}(\Omega)} \le C_l h^{-l},
    \quad l = 0, 1, \ldots, m,
\end{equation}
by Definition \ref{def.MCG} and Assumption A.

The local approximation spaces will satisfy:\smallskip

\noindent {\bf Assumption C.}\ The space $\Psi_j$ contains all
(restrictions to $\omega_j$ of) first order polynomial functions
and there exists  $A > 0$ such that
\begin{equation}\label{eq.2.11}
    \|w\|_{H^l(\omega_j)} \le A \|w\|_{H^l(\omega^*_j)}
\end{equation}
for any $j$, any $w \in \Psi_j$, and any $0 \le l \le m$.
\smallskip

Finally, our last assumptions is the following ``inverse
assumption:''
\smallskip

\noindent {\bf Assumption D.}\ There exists  a constant $B > 0$
\begin{equation}\label{eq.2.13}
    \|w\|_{H^t(\omega_j)} \le B
    d_j^{t-s}\|w\|_{H^s(\omega_j)},
\end{equation}
for any $j$, any $w \in \Psi_j$, and any $0 \le s \le t \le m$.
\smallskip

Assumptions C and D are satisfied, for example, if $\Psi_j$ is the
space of polynomials of degree $\le k-1$, for some fixed $k \ge
2$, see Section \ref{sec.Ver}.

If $A \subset \Omega$ is an open subset, then $A$ satisfies
Assumptions B, C, and D, but not A, in general. We shall hence
need to single out the class of open subsets of $\Omega$
satisfying Assumption A.

\begin{definition}\label{def.adm}\
An open subset $A \subset \Omega$ is called {\em admissible} if
$\sum_{j \in J(A)} \phi_j = 1$ on $A$, where $J(A)$ is the set of
those indices $j$ such that $\omega^*_j \subset A$.
\end{definition}

\begin{remark}\label{rem.admiss}\
It is easy to see that $A \subset \Omega$ is an admissible open
subset if, and only if, $A$ satisfies Assumption A. In particular,
$\Omega$ is an admissible subset of itself (see also Equation
\eqref{eq.res.ABCD}). Moreover, all our results on the set
$\Omega$ extend without change (including the constants) to any
admissible open subset $\Omega_1 \subset \Omega$.
\end{remark}

It is proved in Section \ref{sec.Ver} that enough admissible open
sets exist if $h$ is small enough.  Almost all open sets below
will be admissible. A typical example of an admissible set is
obtained as follows. Fix a subset $J$ of indices $j$ and let $G$
be the set of points where $\sum_{j \in J} \phi_j = 1$. Then the
interior of $G$ is an admissible open subset.

\subsection{Discrete solution\label{ssec.DS}}
Consider the usual bilinear form
\begin{equation}
    B(u, v) := \int_{\Omega} \nabla u \cdot \nabla v dx,
\end{equation}
defined, for example, for $u, v \in H\sp{1}(\Omega)$. Assume, for the
purpose of this discussion, that $\Omega$ is smooth, so that the
boundary value problem \eqref{eq.BVP} makes sense. Then the
solution $u$ of Equation \eqref{eq.BVP} satisfies
\begin{equation}
    B(u, v) = \pairing{g}{v\vert_{\pa \Omega}},
\end{equation}
for $v$ smooth enough. We define then the {\em GFEM--solution} of
Equation \eqref{eq.BVP} to be $u_S \in S$ such that
\begin{equation}
    B(u_S, v_S) = \pairing{g}{v_S\vert_{\pa \Omega}},
\end{equation}
and $\int_{\Omega} u_S(x) dx = 0$. This is possible since the form
$B$ is non-degenerate on the subspace $S_0 \subset S$ consisting
of functions with zero integral over $\Omega$ (recall that we are assuming,
for simplicity, that $\Omega$ is connected). See also Lemma
\ref{lemma.disc}. This also shows that we need the discretization
space $S$ to consist of functions that are smooth enough, which is
the main reason why we are using the Generalized Finite Element
Method in this paper.

\section{Interior estimates for the GFEM\label{sec.IE}}

We continue to assume that $\Omega$ is bounded and connected.
Also, we do {\em not} assume that $\Omega$ has a smooth boundary.
We assume, however, that there are given $\Psi_j$, $\omega_j^*
\subset \omega_j$, and $\phi_j$ be as in the previous section. In
particular, they are assumed to satisfy Assumptions A, B, C, and
D. This will rule out some sets $\Omega$, but any $\Omega$ with
piecewise smooth boundary enjoys this property, see Section
\ref{sec.Ver}. Let $S = S_{GFEM}$ be the resulting Generalized
Finite Element Space.

All the parameters appearing in Assumptions A--D, except $h$, will
be fixed in what follows (\ie $A$, $B$, $C_j$, $\kappa$, $m$, and
$\sigma$ will be fixed). In particular, when we shall say that
certain other constants are ``independent of our choice of $h$ and
GFEM-space $S$,'' we shall understand that all possible choices of
a GFEM--space are allowed, as long as the Assumptions A--D are
satisfied for an arbitrary, but fixed, choice of the parameters
$A$, $B$, $C_j$, $\kappa$, and $\sigma$. In particular, our
``constants'' are allowed to depend on $A$, $B$, $C_j$, $\kappa$,
$m$, and $\sigma$. The parameter $h$ that measures the degree of
the refinement of our covering is allowed, however, to become as
small as we want. Also, all constants will be independent on $N$,
the number of open sets in our covering $\{\omega_j\}$.

All our results below remain true (with the same constants) if we
replace $\Omega$ by an arbitrary admissible open subset $\Omega_1
\subset \Omega$. (Recall that admissible open subsets were
introduced in Definition \ref{def.adm}.)

\subsection{$H^k$-approximation}
We shall need a basic result on the approximation of functions in
$H^{k}(\Omega)$ with elements in the GFEM--space $S$, extending
Theorem \ref{thm.2.1}. Only the case $k =1$ will be needed in this
section, but later on we shall need the general case.

First, let us recall the following standard lemma.

\begin{lemma}\label{lemma.f.cov}\
Let $\psi_j$ be measurable functions defined on $\Omega$. Assume
that there exists an integer $\kappa$ such that a point $x \in
\Omega$ can belong to no more than $\kappa$ of the sets
$\supp(\psi_j)$. Let $f = \sum_j \psi_j$. Then there exists a
constant $C > 0$, depending only on $\kappa$, such that
$\|f\|_{H^l(\Omega)}^2 \le C \sum_{j} \|\psi_j\|_{H^l(\Omega)}^2$.
\end{lemma}

\begin{proof}\ The inequality
\begin{equation}
    |a_1 + a_2 + \ldots + a_M|^2 \le M
    \big( |a_1|^2 + |a_2|^2 + \ldots + |a_M|^2 \big)
\end{equation}
gives the desired result.
\end{proof}

Let $k \in \ZZ_+$. We shall denote as usual
\begin{equation*}
    |u|_{W^{k, \infty}(\Omega)} := \max_{|\alpha| =  k}
    \|\pa^\alpha u\|_{L^{\infty}(\Omega)}, \quad
    \|u\|_{W^{k, \infty}(\Omega)} := \max_{|\alpha| \le k}
    \|\pa^\alpha u\|_{L^{\infty}(\Omega)},
\end{equation*}
$W^{k, \infty}(\Omega) := \{ u, \|u\|_{W^{k, \infty}(\Omega)} <
\infty \}$, and $\|\nabla \omega\|_{W^{k, \infty}(\Omega)} :=
\sum_j \|\pa_j \omega\|_{W^{k, \infty}(\Omega)}$. In particular,
$|u|_{W^{0, \infty}(\Omega)} = \|u\|_{W^{0, \infty}(\Omega)} =
\|u\|_{L^{\infty}(\Omega)}$. We are ready now to prove the
following theorem.

\begin{theorem}\label{thm.approx}\
Assume that, for each $j$, the local approximation spaces $\Psi_j$
contain (the restriction to $\omega$) of the degree $l-1$
polynomials. Then, for any $v \in H^{l}(\Omega)$ and any any $0
\le k \le l$, there exists $w \in S$ such that
\begin{equation*}
    \|v - w\|_{H^{k}(\Omega)} \le Ch^{l-k} \|v\|_{H^{l}(\Omega)}
\end{equation*}
for a constant $C$ independent of our choice of $h$, $S$, and $v
\in S$.
\end{theorem}

Let us notice that, by taking $k = l$ in the above theorem, we
immediate obtain that, using the same notation,
\begin{equation}\label{eq.k=l}
    \|w\|_{H^{l}(\Omega)} \le C \|v\|_{H^{l}(\Omega)}.
\end{equation}

\begin{proof}\ We shall use the notation and the results
from \cite{BrennerScott}[Chapter 4]. Let $Q_{y, f}(x)$ for any $y
\in \omega^*_j$ be the Taylor polynomial of degree $l-1$ at $y$
associated to a smooth function $f$ . Let then
\begin{equation*}
    w_j(f)(x) := (vol(\omega_j^*))^{-1} \int_{\omega^*_j} Q_{y,
    f}(x)dy \in \Psi_j
\end{equation*}
be the  Taylor polynomial of degree $l-1$ averaged over
$\omega^*_j$. (In the terminology of \cite{BrennerScott}, this is
the ``Taylor polynomial of order $l$ averaged over $\omega^*_j$.)
The definition of $w_j(f)$ extends to any $f \in H^{l-1}(\Omega)$.
Then we have the well known Bramble--Hilbert Lemma
\cite{BrennerScott}[Lemma 4.3.8]
\begin{equation}\label{eq.BH}
    |f - w_j(f)|_{H^{k}(\omega_j)} \le Ch^{l-k}
    |f|_{H^{l}(\omega_j)},
\end{equation}
with a constant $C$ depending only on $\sigma$. Let $w_j = w_j(v)
\in \Psi_j$ and $w = \sum_{j = 1}^N \phi_j w_j \in S$. Then, using
also Assumption B and Lemma \ref{lemma.f.cov}, we obtain
\begin{multline*}
    |v - w|_{H\sp{s}(\Omega)}^2 \le C \sum_{j=1}^N
    |\phi_j(v - w_j)|_{H\sp{s}(\omega_j)}^2 \\ \le C
    \sum_{j=1}^N \sum_{i=1}^{s} |\phi_j|_{W^{i, \infty}(\omega_j)}^2
    |(v - w_j)|_{H\sp{s - i}(\omega_j)}^2 \le C \sum_{j=1}^N
    \sum_{i=1}^{s} C_i^2 h^{-2i} h^{2l - 2s +
    2i} |v|_{H\sp{l}(\omega_j)}^2\\
    \le C \kappa h^{2l - 2s} |v|_{H\sp{l}(\Omega)}^2
\end{multline*}
Summing over $0 \le s \le k$ gives the desired result.
\end{proof}

Let us also record, for further use, the following well known
Poincar\'e--Friedrichs inequality \cite{BrennerScott}, Lemma
(4.3.8). (See also \cite{BrennerScott}, Lemma (4.3.14), and
\cite{Ciarlet91}, Equation (2.2), Theorem 14.1, and Theorem 15.3.,
or \cite{Evans, Taylor1}.) The precise statement that we need is
the following.

\begin{theorem}\label{thm.2.2}\
Using the notation of Theorem \ref{thm.2.1}, we have
\begin{equation}\label{eq.2.17}
    \|v\|_{L^2(\omega_j)} \le C_P h |v|_{H^1(\omega_j)},
\end{equation}
for all $v \in H^1(\omega_j)$ satisfying $\int_{\omega_j} v dx =
0$, where $C_P$ depends only on $\sigma$.
\end{theorem}

Theorems \ref{thm.2.1} and \ref{thm.2.2} lead to the estimate
\begin{equation}\label{eq.2.18}
    |u - u_{ap}|_{H^1(\Omega)} \le (2\kappa)^{1/2} (C_1^2C_P^2 +
    C_0^2)\sp{1/2} \big(\sum \|\nabla (u -
    v_j)\|_{L^2(\omega_j)}^2 \big)\sp{1/2},
\end{equation}
provided that $\int_{\omega_j} (u - v_j) dx = 0$.

For $k = 1$, we shall need the following consequence of Theorem
\ref{thm.approx}, which replaces Assumption 9.5 of
\cite{Wahlbin91}. Define
\begin{equation}\label{eq.def.<}
    S^{<} := S_{GFEM} \cap \mathcal{C}_c(\Omega).
\end{equation}
That is, $S^{<}$ consists of the elements of $S_{GFEM}$ with
compact support inside $\Omega$.

\begin{proposition}\label{prop.A9.5}\
Let $\Omega_1 \Subset \Omega$ and $\theta$ be the distance from
$\pa \Omega$ to $\pa \Omega_1$. Then there exists $C > 0$,
independent of $\theta$, $h$, and the GFEM-space $S$, with the
following property. For any $u \in H^2(\Omega)$ with support in
$\Omega_1$, there exists $w \in S^{<}$ such that
\begin{equation*}
    \|u - w\|_{H^1(\Omega)} \le C \theta^{-1} h
    \|u\|_{H^2(\Omega)}.
\end{equation*}
\end{proposition}

\begin{proof}\
Choose $w_j$ and $w$ as in the proof of Theorem \ref{thm.approx}.
If $h < \theta$, then  $w_j = 0$ unless $\omega_j$ intersects
$\Omega_1$, which gives that the closure of $\omega_j$ is
completely contained in $\Omega$. Hence the support of $w$
constructed above is compact. For $h \ge \theta$, it is enough to
choose $C$ large, by Theorem \ref{thm.approx}.
\end{proof}

\subsection{The super-approximation property} Recall the
bilinear form
\begin{equation}
    B(u, v) := \int_{\Omega} \nabla u \cdot \nabla v dx,
    \quad u, v \in H\sp{1}(\Omega).
\end{equation}

Our approach follows the approach from \cite{NitscheSchatz74}, as
presented in \cite{Wahlbin91}[Section 9]. See \cite{BNS,
NitscheSchatz72, SchatzWahlbin77, SchatzWahlbin95} for related
results on approximation in the ``sup''--norm. Recall that $A
\Subset B$ means that the closure of $A$ is a compact set
contained in the interior of $B$. Our main goal in this section is
to prove Theorem \ref{thm.9.2}.

\begin{lemma}\label{lemma.3.2}\
Let $\rho$ be a smooth function on $\omega_j$ and $w \in \Psi_j$. Then
there exists $\tilde w \in \Psi_j$ such that
\begin{equation}
    \|\rho w - \tilde w\|_{H^1(\omega_j)}
    \le \hat C h \|\rho\|_{W\sp{2, \infty}(\omega_j)}
    \|w\|_{H^1(\omega_j)},
\end{equation}
where $\hat C > 0$ may depend only on the dimension $n$ (in
particular, it is independent of $w$, $\rho$, or $j$).
\end{lemma}

\begin{proof}\
We shall use the inner product induced from $H^1(\omega_j)$. Let
$\rho \in W\sp{2, \infty}(\omega_j)$ be given.

To prove the lemma, we shall assume first that $w \in \Psi_j$ is a
constant. Let $l$ be the degree one Taylor polynomial
approximation of $\rho$ at the center of the ball $\omega_j\sp*$.
Then $l \in \Psi_j$, because first order polynomials are in
$\Psi_j$ (Assumption C) and we have
\begin{equation*}
    \|\rho - l\|_{W\sp{1, \infty}(\omega_j)}
    \le \hat C h\|\rho\|_{W\sp{2, \infty}(\omega_j)},
\end{equation*}
with $\hat C > 0$ a constant depending only on the dimension $n$.
(This is where the condition $h \le 1$ is used.) Choose $\tilde w
= l w$. Then
\begin{multline*}
    \|\rho w - \tilde w\|_{H^1(\omega_j)} =
    \|\rho w - l w\|_{H^1(\omega_j)} \le
    \|\rho - l\|_{W\sp{1, \infty}(\omega_j)}
    \|w\|_{H^1(\omega_j)} \\ \le \hat C h
    \|\rho\|_{W\sp{2, \infty}(\omega_j)}
    \|w\|_{H^1(\omega_j)}.
\end{multline*}

Assume now that $w \in \Psi_j$ is such that $\innerp{w}{1} = 0$,
that is, $w$ is orthogonal in $H^1(\omega_j)$ to the subspace
generated by constants. We then write
\begin{equation*}
    \rho = \tilde \rho + \rho^*,
\end{equation*}
where $\tilde \rho$ is a constant function (say the value of
$\rho$ at the center of $\omega\sp*$) and
\begin{equation}\label{eq.est}
    \|\rho^*\|_{L^\infty(\omega_j)} \le \hat C h
    \|\nabla \rho\|_{L^\infty(\omega_j)},
\end{equation}
where $\hat C$ is a constant depending only on the dimension $n$.
We shall choose then $\tilde w = \tilde \rho w \in \Psi_j$, which
makes sense since $\Psi_j$ is a vector space. Then
\begin{multline*}
    \|\rho w - \tilde w\|_{H\sp1(\omega_j)} =
    \|\rho^* w\|_{H\sp1(\omega_j)} \le C
    \|\nabla \rho^*\|_{L\sp{\infty}(\omega_j)}
    \|w\|_{L^2(\omega_j)} \\ + C
    \|\rho\sp*\|_{L\sp{\infty}(\omega_j)}
    \|w\|_{H^1(\omega_j)} \le \hat C h
    \|\nabla \rho\|_{L^\infty(\omega_j)}
    \|w\|_{H^1(\omega_j)},
\end{multline*}
where in the last step we have used the Poincar\'e--Friedrichs
inequality for $\omega_j$ (Theorem \ref{thm.2.2}) to estimate
$\|u\|_{L^2(\omega_j)}$ and Equation \eqref{eq.est} above to
estimate the second term. Here $\hat C$ is again a constant that
may depend only on the dimension $n$.

For a general $w \in \Psi_j$, we decompose $w = w_1 + w_2$ with
$w_1$ a constant and $w_2$ orthogonal to the space of constants
and choose $\tilde w_1$ and $\tilde w_2$ as above. Then
\begin{multline*}
    \|\rho w - \tilde w_1 - \tilde w_2\|_{H\sp1(\omega_j)}
    \le  \|\rho w_1 - \tilde w_1\|_{H\sp1(\omega_j)} +
    \|\rho w_2 - \tilde w_2\|_{H\sp1(\omega_j)} \\ \le
    \hat C h \|\rho\|_{W\sp{2, \infty}(\omega_j)}
    \|w_1\|_{H^1(\omega_j)} + \hat C h
    \|\rho\|_{W\sp{2, \infty}(\omega_j)}
    \|w_2\|_{H^1(\omega_j)} \\ \le
    C h \|\rho\|_{W\sp{2, \infty}(\omega_j)}
    \|w\|_{H^1(\omega_j)}.
\end{multline*}
The lemma is now proved.
\end{proof}

An important technical step in our proof of the Theorem
\ref{thm.9.2} is the following ``super-approximation'' result.

\begin{proposition}\label{prop.lemma.3.1}\
Let $\rho \in W\sp{2, \infty}(\Omega)$ and $w \in S = S_{GFEM}$.
Then there exists $\tilde w \in S$ such that
\begin{equation}\label{eq.3.2}
    \|\rho w - \tilde w\|_{H\sp1(\Omega)} \le C h
    \|\rho\|_{W\sp{2, \infty}(\Omega)} \|w\|_{H\sp1(\Omega)},
\end{equation}
where $C$ is independent of $h$ and our choice of GFEM--space $S$.
If $u$ has  support in $\Omega_1 \Subset \Omega$, then we can
chose $\tilde w$ to have compact support in $\Omega$ and $C \le C'
\theta^{-1}$, where $\theta$ is the distance from $\pa \Omega$ to
$\pa \Omega_1$ and $C'$ is independent of $\theta$, $h$, and $S$.
\end{proposition}

As explained above, the constant $C$ may depend, however, on the
parameters $A$, $B$, $C_j$, $\kappa$, $m$, $\sigma$, and $\theta$,
but is independent of $h$ or of the number $N$ of sets
$\{\omega_j\}$.

\begin{proof}\ Let
\begin{equation}\label{eq.3.7}
    w = \sum_{j = 1}^N \phi_j w_j \in S = S_{GFEM}, \quad
    w_j \in \Psi_j.
\end{equation}
Let $\tilde w_j$ be the orthogonal projection of $\rho w_j$ onto
$\Psi_j$ in the inner product of $H\sp{1}(\omega_j)$. Lemma
\ref{lemma.3.2} then shows that
\begin{equation}
    \|\rho w_j - \tilde w_j \|_{H^1(\omega_j)}
    \le \hat C h \|\rho\|_{W\sp{2, \infty}(\Omega)}
    \|w_j\|_{H\sp1(\omega_j)}.
\end{equation}
Moreover, we have that $\int_{\omega_j}(\rho w_j - \tilde w_j)dx =
0$ because the constant functions are in $\Psi_j$ and $\rho w_j -
\tilde w_j$ is orthogonal to $\Psi_j$.

Let $\tilde w := \sum_{j = 1}^N \phi_j \tilde w_j$. Then $\|\nabla
\phi_j\|_{L\sp\infty(\omega_j)} \le C_1/h$ by Equation
\eqref{eq.2.10} and
\begin{equation*}
    \|\rho w_j - \tilde w_j\|_{L\sp2(\omega_j)} \le
    C_P h \|\rho w_j  - \tilde w_j\|_{H\sp1(\omega_j)} \le
    C_P \hat C h^{2} \|\rho\|_{W\sp{2, \infty}(\Omega)}
    \|w_j\|_{H\sp1(\omega_j)},
\end{equation*}
by the Poincar\'e-Friedrichs inequality (Theorem \ref{thm.2.2}),
and hence
\begin{multline*}
    \|\rho w - \tilde w\|_{H\sp{1}(\Omega)}^2  =
    \|\sum_{j = 1}^N \phi_j (\rho w_j - \tilde
    w_j)\|_{H\sp{1}(\Omega)}^2 \le C \sum_{j = 1}^N \Big (
    \|\phi_j\|_{L^\infty(\omega_j)}^2
    \|\rho w_j - \tilde w_j\|_{H\sp{1}(\omega_j)}^2 \\
    + \|\nabla \phi_j\|_{L\sp\infty(\omega_j)}^2
    \|\rho w_j - \tilde w_j\|_{L\sp2(\omega_j)}^2 \Big)
    \le C h ^2 \|\rho\|_{W\sp{2, \infty}(\Omega)}^2
    \sum_{j = 1}^N \|w_j\|_{H\sp1(\omega_j)}^2 ,
\end{multline*}
where for the first inequality we have used also Lemma
\ref{lemma.f.cov}. The result will follow now if we can prove that
$\sum_{j = 1}^N \|w_j\|_{H\sp1(\omega_j)}^2  \le C
\|w\|_{H\sp1(\Omega)}^2$, for any $w = \sum_{j = 1}^N \phi_j w_j$,
as above and $C$ a constant independent of $h$ and $S$. Indeed, we
have
\begin{equation*}
    \|w\|_{H\sp1(\Omega)}^2 \ge \sum_{j = 1}^N
    \|w_j\|_{H\sp1(\omega_j^*)}^2 \ge A^{2} \sum_{j = 1}^N
    \|w_j\|_{H\sp1(\omega_j)}^2,
\end{equation*}
by Assumption C ($A$ is the constant appearing in that
assumption).

The proof of the last part is completed as in Proposition
\ref{prop.A9.5}.
\end{proof}

\subsection{Estimates on ``discrete--harmonic'' functions}
We shall also need the following ``inverse property,'' which is
\CORR somewhat similar to Assumption A.3. in
\cite{NitscheSchatz74} or Assumption 9.2 in \cite{Wahlbin91}.

\begin{lemma}\label{lemma.4.1}\
There exists $C > 0$, independent of $h$, $u$, and the GFEM--space
$S$ such that
\begin{equation}
    \|w\|_{H\sp{j}(\Omega)} \le C h^{i - j} \|w\|_{H\sp{i}(\Omega)},
\end{equation}
for all $0 \le i \le j \le m$.
\end{lemma}

Recall that the constant $m$ is the fixed integer appearing in
Assumptions A--D.

\begin{proof}\
If $i = j$, we can take $C = 1$.  Let $w = \sum_{k = 1}^N \phi_k
w_k$, with $w_k \in \Psi_k$. Then Lemma \ref{lemma.f.cov} and
Assumptions A--D give
\begin{multline*}
    \|w\|_{H\sp{j}(\Omega)}^2 \le \sum_{k = 1}^N \sum_{l = 0}^j
    \|\phi_k\|_{W^{l, \infty}(\omega_k)}^2 \|w_k\|_{H\sp{j -
    l}(\omega_k)}^2 \\ \le \sum_{k = 1}^N\sum_{l = 0}^j C_l^2
    h^{-2l} B^2 h^{2l-2j} \|w_k\|_{L^2(\omega_k)}^2 \le C^2
    A^2 h^{-2j} \sum_{k = 1}^N \|w_k\|_{L^2(\omega_k^*)}^2 \\ \le
    C^2 h^{-2j} \|w\|_{L^2(\Omega)}^2.
\end{multline*}
This proves the result for $i = 0$. For an arbitrary $0 \le i \le
j$, the result follows by interpolation.
\end{proof}

\CORR The rest of this section follows closely the approach in the
paper of Nitsche and Schatz \cite{NitscheSchatz74}, relying also
from the survey paper \cite{Wahlbin91} (which in turn is based on
the paper by Nitsche and Schatz). There are, however, some
differences in the assumptions that we are using, so we include
complete proofs for the convenience of the reader. For instance,
\CORR the following corollary of Lemma \ref{lemma.4.1} plays the
role of Assumption A.3. in the Nitsche--Schatz article
\cite{NitscheSchatz74}, respectively, of the Assumption 9.2
(Inverse assumption) in Wahlbin's article. Also, the following
lemma is an analog of Lemma 5.2 of \cite{NitscheSchatz74},
respectively, of Lemma 9.1 of \cite{Wahlbin91}. Recall that all
the above results remain true if $\Omega$ is replaced by an
admissible open subset $\Omega_1 \subset \Omega$.

\begin{corollary}\label{cor.inv}\
Let $w \in S$. Then there exists a constant $C > 0$ such that
\begin{equation*}
    \|w\|_{H\sp{-i}(\Omega)} \le C h ^{i-j}
    \|w\|_{H\sp{-j}(\Omega)},
\end{equation*}
for any $0 \le i \le j \le m$.
\end{corollary}

\begin{proof}\ If $i = j$ we can take $C = 1$. Assume next that $i
=0$. Then
\begin{equation*}
    \|w\|_{H\sp{-j}(\Omega)} \! = \sup \frac{|(w,
    \phi)|}{\;\|\phi\|_{H\sp{j}(\Omega)}}
    \ge \frac{(w, w)}{\;\|w\|_{H\sp{j}(\Omega)}} =
    \frac{\|w\|^2_{L\sp{2}(\Omega)}}{\|w\|_{H\sp{j}(\Omega)}}
    \ge C^{-1}h^{j} \|w\|_{L\sp{2}(\Omega)},
\end{equation*}
by Lemma \ref{lemma.4.1}. For the other values of $j$, the result
follows by interpolation.
\end{proof}

We shall denote $A \Subset B$ if $\overline{A}$, the closure of $A$ in
$\RR^2$, is a compact subset of the interior of $B$.  Also, we shall
denote
\begin{equation}\label{eq.def.S<}
    S^{<}(A) := \{ u \in S_{GFEM}, \supp(u) \Subset A\}
\end{equation}
for any {\em admissible} open subset $A \subset \Omega$. In
particular, $S^{<}(\Omega) = S^{<}$.

\begin{lemma}\label{lemma.9.1}\
Let $A \Subset A_1 \Subset \Omega$ be admissible open sets and
$\theta = dist(\pa A, \pa A_1)$.  Then there exists $C > 0$,
independent of $h$, $S$, $A$, and $A_1$ with the following
property. If $w \in S = S_{GFEM}$ and
\begin{equation}\label{eq.9.10}
    B(w, \chi) = 0, \quad \text{for all }\chi \in S^{<}(A_1),
\end{equation}
then, for $h$ small enough, $\|w\|_{H^1(A)} \le C
\|w\|_{L^2(A_1)}$, with $C$ depending only on $\theta$ and not on
$A$, $A_1$, $h$, or $S$.
\end{lemma}

\begin{proof}\
The proof is the same as the one in \cite{NitscheSchatz74,
Wahlbin91}, using Lemma \ref{lemma.4.1} in place of Assumption
A.3, respectively Assumption 9.2 (``Inverse assumption''), and
Proposition \ref{prop.lemma.3.1} in place of Assumption A.2,
respectively Assumption 9.1 (``Superapproximation''), of
\cite{NitscheSchatz74}, respectively \cite{Wahlbin91}. The
constants ``$C$'' below are allowed to depend on $\theta$.

Let us chose $A \Subset A_0 \Subset A_1$ admissible open sets such
that the distances between the boundaries of these sets are $\ge
\theta/C$. This is possible if $h < c\theta$ (see Section
\ref{sec.Ver}). Also, let $\omega \in \CIc(A_0)$, with $\omega =
1$ on $A$, $\omega \ge 0$ on $A_0$, and $\|\omega\|_{W^{k,
\infty}(\Omega)} \le C \theta^{-k}$, for $k = 0, 1$. Then, by
Equation \eqref{eq.9.10}, we obtain
\begin{multline*}
    \|\nabla w\|_{L^2(A)}^2
    \le (\nabla w, \omega \nabla w)
    = (\nabla w, \nabla (\omega w))
    - (\nabla w, \nabla (\omega) w) \\
    = (\nabla w, \nabla (\omega w - \psi))
    + \ha (w, (\Delta \omega) w),
\end{multline*}
where the inner products are in $L^2(A_0)$ and $\psi \in
S^{<}(A_0)$. Proposition \ref{prop.lemma.3.1} then gives
\begin{equation*}
    \|\nabla w\|_{L^2(A)}^2 \le C h
    \|\omega\|_{W^{2, \infty}(\Omega)}
    \|w\|_{H^1(A_0)}^2 + C \|w\|_{L^2(A_0)}^2,
\end{equation*}
which, in turn, implies
\begin{equation}\label{eq.eq1}
    \|w\|_{H\sp1(A)} \le C \big ( h^{1/2}
    \|w\|_{H^1(A_0)} + \|w\|_{L^2(A_0)}).
\end{equation}
We now repeat the argument for $A_0 \Subset A_1 \Subset \Omega$ (and
$A$ replaced by $A_0$ and $A_0$ replaced by $A_1$), which gives
\begin{equation}\label{eq.eq2}
    \|w\|_{H^1(A_0)} \le C \big( h^{1/2}
    \|w\|_{H^1(A_1)} + \|w\|_{L^2(A_1)} \big).
\end{equation}
Combining Equations \eqref{eq.eq1} and \eqref{eq.eq2} and using
also $h \le 1$, we obtain
\begin{equation}\label{eq.eq3}
    \|w\|_{H\sp1(A)} \le C \big(h \|w\|_{H^1(A_1)}
    + C \|w\|_{L^2(A_1)}\big).
\end{equation}
Since $A_1$ is admissible, we can use Lemma \ref{lemma.4.1} with
$\Omega$ replaced with $A_1$ to obtain $h \|w\|_{H^1(A_1)} \le C
\|w\|_{L^2(A_1)}$, with $C$ independent of $h$, $w$, and the GFEM
space $S$ (as long as $S$ satisfies the Assumptions A--D). Then
\begin{equation}\label{eq.eq4}
    \|w\|_{H\sp1(A)} \le C \|w\|_{L^2(A_1)}.
\end{equation}
The proof is now complete.
\end{proof}

We shall need the following simple estimate.

\begin{lemma}\label{lemma.conv.N}\
Let $\Phi(x) = \log |x|$ if $n = 2$, $\Phi(x) = |x|^{2-n}$, if $n
\neq 2$. Let $U$ be a fixed bounded open subset of $\RR^n$. Then
there exists $C > 0$, which depends only on $U$, such that
\begin{equation*}
    \Phi*u(x) := \int_{y} \Phi(x-y) u(y) dy
\end{equation*}
satisfies
\begin{equation*}
    \|\Phi*u\|_{H^{l+2}(U)} \le C \|u\|_{H^{l}(U)},
\end{equation*}
for any $l \in \RR$ and any $u \in \CIc(U)$.
\end{lemma}

\begin{proof}\ Let $\omega \in \CIc(\RR^n)$ be equal to $1$
on $U$. Then $Tu(x) = \omega(x) [\Phi*(\omega u)](x)$ is a
pseudodifferential operator of order $-2$ with compactly supported
distribution kernel. Hence it is bounded as a map $H^{l}(\RR^n)
\to H^{l+2}(\RR^n)$ \cite{hor3, Taylor2}. The statement follows by
restricting $T$ to $\CIc(U) \subset H^{l}(\RR^n)$.
\end{proof}

We define
\begin{equation}
    \|u\|_{H_0^{-s}(U)} = \sup \frac{|(u, v)|}{\|v\|_{H\sp{l}(U)}}
    \le \|u\|_{H^{-s}(U)}, \quad 0 \neq v \in \CIc(A)
\end{equation}
for any open set $U$, any $u \in L^2(U)$, and any $s > 0$. We
define $H_0^{-s}(U)$ to be the completion of $L^2(U)$ in the norm
$\|u\|_{H_0^{-s}(U)}$. Then $H_0^{-s}(U)$, $s
> 0$, identifies with the dual of $H_0^{s}(U)$.

We now prove the following lemma.

\begin{lemma}\label{lemma.9.2}\
We keep the notation and assumption of Lemma \ref{lemma.9.1}. In
particular, we assume that $w \in S$ satisfies Equation
\eqref{eq.9.10}. Then, for $h$ small enough,
\begin{equation}
    \|w\|_{L^2(A)} \le C \|w\|_{H^{-m}(A_1)},
\end{equation}
where $C$ is a constant depending only on the distance $\theta$
from $\pa A$ to $\pa A_1$ and not on $A$, $A_1$, $h$, or $S$.
\end{lemma}

Combining Lemmata \ref{lemma.9.1} and \ref{lemma.9.2}, we obtain
\begin{equation}\label{eq.ref}
    \|w\|_{H\sp{1}(A)} \le C \|w\|_{H^{-m}(A_1)},
\end{equation}
for $h$ small enough and any $w \in S$ satisfying the assumptions
of Lemma \ref{lemma.9.1}.

\begin{proof}\
Let $A \Subset B_0 \Subset B_1 \Subset A_1$ be such that the
distances between the boundaries of these sets are $\ge \theta/C$.
This is possible if $h < c\theta$.  Note that by Lemma
\ref{lemma.9.1}, we have
\begin{equation}\label{eq.lemma.9.1}
    \|w\|_{H\sp{1}(B_1)} \le C \|w\|_{L^2(A_1)}.
\end{equation}

For any $v \in \CIc(A)$, let $V := c_n \Phi *v \in H^{l + 2}(B_1)$,
where $c_n$ is chosen such that $\Delta V = v$ (see
\cite{Evans}). Lemma \ref{lemma.conv.N} then gives
\begin{equation}\label{eq.as.9.4}
    \|V\|_{H\sp{l + 2}(B_1)} \le C \|v\|_{H\sp{l}(B_1)} = C
    \|v\|_{H\sp{l}(A)}, \quad l \in \ZZ_+,
\end{equation}
for some constant $C$ that depends only on $B_1$. Let $\omega \in
\CIc(B_0)$ with $\omega = 1$ on $A$ and $\|\omega\|_{W^{k,
\infty}(\Omega)} \le C \theta^{-k}$, for $0 \le k \le m$. Since
$\omega V \in \CIc(B_0)$, we know from Proposition \ref{prop.A9.5}
that there exists $\chi \in S^{<}(B_1)$ such that
\begin{equation}\label{eq.as.9.5}
    \|\omega V - \chi\|_{H\sp{1}(B_1)} \le C h \|\omega
    V\|_{H\sp{2}(B_1)} \le C h \| V \|_{H\sp{2}(B_1)} \le C h
    \| v \|_{L^2(A)}.
\end{equation}

Then, for any $v \in \CIc(A)$,
\begin{multline*}
    (w, v)_{A} = (\omega w, v)_{A} = \int_{A} \omega w
    \Delta V dx = \int_{B_0} \omega w \Delta V dx =
    \int_{B_0}\nabla (\omega w)\cdot \nabla V dx \\ = \int_{B_0}
    w \big(2 \nabla \omega \cdot \nabla V - V\Delta \omega \big)
    dx + \int_{B_1}\nabla w \cdot \nabla (\omega V - \chi) dx,
\end{multline*}
for any $\chi \in S^{<}(B_1)$, where the inner products are
calculated on the indicated sets. Then, by combining Equations
\eqref{eq.lemma.9.1}, \eqref{eq.as.9.4}, and \eqref{eq.as.9.5}, as
well as Lemma \ref{lemma.9.1}, we obtain, for all $l \ge 0$,
\begin{multline*}
    |(w, v)_{A}| \le C \|w\|_{H_0\sp{-l -1}(B_0)}
    \|V\|_{H\sp{l + 2}(B_1)} + h \|w\|_{H\sp{1}(B_1)}
    \|V\|_{H\sp{2}(B_1)}\\ \le C \big(\|w\|_{H_0\sp{-l -1}(A_1)}
    + h\|w\|_{L\sp{2}(A_1)} \big) \|v\|_{H\sp{l}(A)},
\end{multline*}
and hence
\begin{equation}\label{eq.to.it}
    \|w\|_{H_0\sp{-l}(A)} \leq  C \big(\|w\|_{H_0\sp{-l
    -1}(A_1)} + h\|w\|_{L\sp{2}(A_1)} \big) .
\end{equation}

Let us now choose a sequence of open sets $A \Subset B_1 \Subset
B_2 \Subset \ldots \Subset B_m \Subset \Omega$ with all distances
between the boundaries greater or equal $c\theta$. Changing
notation and iterating Equation \eqref{eq.to.it}, we obtain
\begin{multline}\label{eq.m.to.it}
    \|w\|_{L\sp{2}(A)} \leq  C
    \big(\|w\|_{H_0\sp{-1}(B_1)}
    + h\|w\|_{L\sp{2}(B_1)} \big) \\ \le C
    \big(\|w\|_{H_0\sp{-2}(B_2)}
    + h\|w\|_{L\sp{2}(B_2)} \big) \le \ldots \\
    \le  C \big(\|w\|_{H_0\sp{-m}(B_m)}
    + h\|w\|_{L\sp{2}(B_m)} \big).
\end{multline}

We now repeat the above reasoning. We change notation again, so
that, this time, $B_m$ becomes $B_1$, then we chose as before a
sequence of open sets
\begin{equation*}
    A \Subset B_1 \Subset B_2 \Subset \ldots \Subset B_m
    \Subset \Omega
\end{equation*}
with all distances between the boundaries
greater or equal $c\theta$. Then we iterate Equation
\eqref{eq.m.to.it}, and obtain,
\begin{multline*}
    \|w\|_{L\sp{2}(A)} \leq  C
    \big(\|w\|_{H_0\sp{-m}(B_1)}
    + h\|w\|_{L\sp{2}(B_1)} \big) \\ \le C
    \big(\|w\|_{H_0\sp{-m}(B_2)}
    + h^2\|w\|_{L\sp{2}(B_2)} \big) \le \ldots
    \le  C \big(\|w\|_{H_0\sp{-m}(B_m)}
    + h^m\|w\|_{L\sp{2}(B_m)} \big) \\
    \le C \big(\|w\|_{H_0\sp{-m}(B_m)}
    + \|w\|_{H\sp{-m}(B_m)} \big) \le C
    \|w\|_{H\sp{-m}(B_m)},
\end{multline*}
where at the end we have used the inverse property
$\|w\|_{L\sp{2}(U)} \le h^{-m}\|w\|_{H\sp{-m}(U)}$ for any
admissible open set $U \subset \Omega$ (see Corollary
\ref{cor.inv}). The proof is now complete.
\end{proof}

\subsection{The interior error estimate}
The following result, the main result of this section, \CORR is an
analog of \cite{NitscheSchatz74}[Theorem 5.1] and of
\cite{Wahlbin91}[Theorem 9.2].

\begin{theorem}\label{thm.9.2}\
Let $A \Subset B \subset \Omega$ be admissible open sets. Then
there exists $C > 0$ with the following property. If $u \in
H^1(\Omega)$ and $u_S \in S$ are such that $B(u - u_S, \chi) = 0$
for all $\chi \in S^{<} := S^{<}(\Omega)$, then for $h$ small
enough,
\begin{equation*}
    \|u - u_S\|_{H^1(A)} \le C \Big( \inf_{\chi \in S}\|u -
    \chi\|_{H^1(B)} + \|u - u_S\|_{H\sp{-m}(B)} \Big).
\end{equation*}
The constant $C$ depends only on the distance $\theta = dist(\pa
A, \pa B)$ and not on $h$, $S$, or the sets $A$ and $B$.
\end{theorem}

\begin{proof}\
Let $A \Subset A_1 \Subset A_2 \Subset B \subset \Omega$.  Choose
$\omega \in \CIc(A_2)$ such that $\omega = 1$ on $A_1$. Let $P_1$
be the $H^1(\Omega)$ orthogonal projection onto $S^{<}(A_1)
\subset S \subset H\sp{1}(\Omega)$. Then on $A_1$
\begin{equation}\label{eq.one}
    u - u_S = \big (\omega u - P_1(\omega u) \big) +
    \big (P_1(\omega u)  - u_S\big).
\end{equation}
Then, by the general properties of orthogonal projections, we have
\begin{equation}\label{eq.two}
    \|\omega u - P_1(\omega u)\|_{H^1(\Omega)}
    \le \|\omega u\|_{H^1(\Omega)} \le
    C \|u\|_{H^1(B)}.
\end{equation}
Hence
\begin{equation}\label{eq.three}
    \|\omega u - P_1(\omega u)\|_{H\sp{-m}(A_1)} \le
    \|\omega u - P_1(\omega u)\|_{H\sp{1}(\Omega)}
    \le C \|u\|_{H^1(B)}.
\end{equation}

Let $w = P_1(\omega u) - u_S$. Then $B(w, \chi) = B(\omega u -
u_S, \chi) = B(u - u_S, \chi) = 0$, for all $\chi \in S^{<}(A_1)$,
and hence $w$ satisfies the assumptions of Lemmata \ref{lemma.9.1}
and \ref{lemma.9.2}. From this, using also Equations
\eqref{eq.ref}, \eqref{eq.one}, and \eqref{eq.three}, we obtain
\begin{multline}\label{eq.four}
    \|w\|_{H\sp{1}(A)} \le C \|w\|_{H\sp{-m}(A_1)}
    \le \|\omega u - P_1(\omega u)\|_{H\sp{-m}(A_1)} +
    \|u - u_S\|_{H\sp{-m}(A_1)}  \\ \le  \|u\|_{H^1(B)} +
    \|u - u_S\|_{H\sp{-m}(A_1)}.
\end{multline}

Equations (\ref{eq.one}--\ref{eq.four}) then give
\begin{eqnarray*}
    \begin{array}{rll}
    \|u - u_S\|_{H\sp{-m}(A)}\! & \le
    \|\omega u - P_1(\omega u)\|_{H\sp{-m}(A)} +
    \|w\|_{H\sp{-m}(A)} & \text{by
    \eqref{eq.one}} \\ & \le C \|u\|_{H\sp{1}(B)} +
    \|u - u_S\|_{H\sp{-m}(A_1)} & \text{by
    \eqref{eq.two}--\eqref{eq.four}.}
    \end{array}
\end{eqnarray*}
The desired result follows by replacing $u$ and $u_S$ with $u -
\chi$ and, respectively, $u_S - \chi$, with $\chi$ in $S =
S_{GFEM}$.
\end{proof}

\section{Approximate solution of the Laplace equation with distribution
boundary conditions using the GFEM\label{sec.AS}}

We shall consider the same setting as in the previous sections,
and we shall further assume that $\Omega$ is a smooth domain. In
particular, $S$ will be the GFEM--space associated to a partition
of unity $\{\phi_j\}$ subordinated to the covering $\{\omega_j\}$
of $\Omega$ and local approximation spaces $\{\Psi_j\}$. We shall
continue to assume that Assumptions A--D are satisfied, for a
fixed choice of constants $A$, $B$, $C_j$, $\kappa$, $m$, and
$\sigma$. In particular, the ``constants'' below are allowed to
depend on these parameters, but are not allowed to depend on $h$
or the specific choice of the GFEM--space $S$, as long as the
constants above remain the same.

We shall denote by $\nu$ the outer unit normal vector to $\pa
\Omega$. By $\pa_\nu u(x)$ we shall denote the directional
derivative of a function $u$ in the direction of $\nu$, at some
point $x$ on the boundary.

Let $u \in H^{1-k}(\Omega)$. We want to make precise in what sense
we shall say that ``$\Delta u = 0$ as a distribution on
$\Omega$.'' Recall first that the space $H^{1-k}(\Omega)$, $k \in
\ZZ_+$, $k > 1$, was defined as the dual of $H\sp{k-1}(\Omega)$
(see Section \ref{sec.prel}). We shall write $\pairing{u}{v} :=
u(v) \in \CC$ for any $u \in H^{1-k}(\Omega)$ and $v \in
H^{k-1}(\Omega)$ for the ``value of $u$ evaluated at $v$.'' We can
hence define by duality $\pa_j := - \pa_j^* : H^{1-k}(\Omega) \to
H\sp{-k}(\Omega)$. This leads to a definition of $\Delta u \in
H^{-1 -k}(\Omega)$, for any $u \in H^{1 -k}(\Omega)$ by
\begin{equation*}
    \pairing{\Delta u}{v} := \pairing{u}{\Delta v},\quad \text{ for
    any } v \in H^{1-k}(\Omega).
\end{equation*}
However, this turns out to be too strong a condition. Instead, we
shall require
\begin{equation}\label{eq.0dist}
    \pairing{\Delta u}{\phi} := \pairing{u}{\Delta \phi} = 0,
    \quad \text{for any } \phi \in \CIc(\Omega).
\end{equation}

We shall say that {\em $\Delta u = 0$ as a distribution on
$\Omega$} whenever Equation \eqref{eq.0dist} is satisfied.

Typically, $u$ as above will arise as a solution of a boundary
value problem, for example, as a solution of the boundary value
problem \eqref{eq.BVP}. In \cite{hor63, LionsMagenes1,
StrichartzHp} it was explained how to define the traces (or
restrictions) $u \vert_{\pa \Omega}$ and $\pa_\nu u\vert_{\pa
\Omega}$ for any $u \in H^{1-k}(\Omega)$ satisfying $\Delta u = 0$
as a distribution on $\Omega$. More generally, we define $u
\vert_{\pa \Omega} \in H\sp{1/2 -k}(\Omega)$ and $\pa_\nu
u\vert_{\pa \Omega} \in H\sp{-1/2 -k}(\Omega)$ by linearity, for
$u = u_1 + u_2$, where $u_1 \in H\sp{1-k}(\Omega)$, $\Delta u_1 =
0$ as a distribution on $\Omega$, and $u_2 \in H\sp{1}(\Omega)$.
We use this to define $B(u, v)$ by
\begin{equation}\label{eq.def.BUv}
    B(u, v) := - \pairing{u}{\Delta v} +
    \pairing{ u\vert_{\pa \Omega}}{\pa_\nu v},
    \quad \text{for any } v \in H^{1+k}(\Omega),
\end{equation}
for $u = u_1 + u_2$ as above. In view of Green's formula (see
\cite{Evans}, for example), $B(u, v) = \int_{\Omega} \nabla u
\cdot \nabla v dx$ if $u, v \in H^1(\Omega)$, as originally
defined.  It is not clear how to define $B(u, v)$ for arbitrary $u
\in H\sp{1-k}(\Omega)$, since $\pairing{u}{\Delta v}$ is defined
but the traces of $u$ may not defined in general.

{\em From now on, we shall fix $u$ such that, and
\begin{equation}\label{eq.assump.u}
    \Delta u = 0, \quad u \in H\sp{1-k}(\Omega),
\end{equation}
as a distribution on $\Omega$, where $k \in \ZZ_+$, $m +1 \ge k
> 0$, is also fixed}. We do not assume that $\Delta u = 0$ in
$H\sp{-1-k}(\Omega)$ (\ie we do not assume $\pairing{u}{\Delta v}
= 0$ for all $v \in H\sp{1+k}(\Omega)$, we only assume
$\pairing{u}{\Delta v} = 0$ for all $\phi \in \CIc(\Omega)$). We
shall also assume that
\begin{equation}\label{eq.assump.u2}
    \pairing{u}{1} = \pairing{\pa_\nu u\vert_{\pa \Omega}}{1} = 0.
\end{equation}

We have the following.

\begin{lemma}\label{lemma.disc}\
There exists a unique $u_S \in S$ such that $\pairing{u_S}{1} = 0$
and
\begin{equation}\label{eq.disc}
    B(u - u_S, v_S) = 0,
\end{equation}
for all $v_S \in S$.
\end{lemma}

\begin{proof}\
Let $S_0$ be the subspace of the GFEM-space $S$ consisting of
functions $\chi_0 \in S$ with $\int_{\Omega} \chi(x)dx = 0$. The
bilinear form $B$ is non-degenerate on $S_0$. This gives the
existence of a unique $u_S \in S_0$ such that Equation
\eqref{eq.disc} is satisfied for all $v_S \in S_0$. Since $S = S_0
+ \CC$, the result follows from $B(u, 1) = B(u_S, 1) = 0$.
\end{proof}

We also have the following simple estimate.

\begin{lemma}\label{lemma.cont.Bu}\ With $u$ as in Equation
\eqref{eq.assump.u} above, we have
\begin{equation*}
    |B(u, v)| \le C\|u\|_{H\sp{1-k}(\Omega)}
    \|v\|_{H\sp{1+k}(\Omega)},
\end{equation*}
for any $v \in H\sp{1-k}(\Omega)$ and a constant $C$ depending
only on $\Omega$.
\end{lemma}

\begin{proof}\
By definition, we have
\begin{multline}\label{cont.Bu}
    |B(u, v)| = |- \pairing{u}{\Delta v} + \pairing{u\vert_{\pa
    \Omega}}{\pa_\nu v}| \le |\pairing{u}{\Delta v}| +
    |\pairing{u\vert_{\pa \Omega}}{\pa_\nu v}| \\ \le
    \|u\|_{H\sp{1-k}(\Omega)} \|\Delta v\|_{H\sp{-1+k}(\Omega)} +
    \|u\|_{H\sp{1/2-k}(\pa \Omega)} \|v\|_{H\sp{-1/2+k}(\pa
    \Omega)}\\ \le C \|u\|_{H\sp{1-k}(\Omega)}
    \|v\|_{H\sp{1+k}(\Omega)}.
\end{multline}
This completes the proof.
\end{proof}

We continue with more lemmata. Recall that $k \in \ZZ_+$.

\begin{lemma}\label{lemma.2}\ We have
$\|u_S\|_{H\sp{1}(\Omega)} \le C h^{-k} \|u\|_{H\sp{1-k}(\Omega)}$
for a constant $C$ depending only on $\Omega$.
\end{lemma}

\begin{proof}\ The Poincar\'e-Friedrichs inequality and Lemma
\ref{lemma.cont.Bu} give
\begin{multline*}
    \|u_S\|_{H\sp{1}(\Omega)}^2 \le C B(u_S, u_S)
    = C B(u, u_S)  \le C \|u\|_{H\sp{1-k}(\Omega)}
    \|u_S\|_{H\sp{1+k}(\Omega)} \\ \le C h^{-k}
    \|u\|_{H\sp{1-k}(\Omega)} \|u_S\|_{H\sp{1}(\Omega)},
\end{multline*}
where in the last inequality we have used Lemma \ref{lemma.4.1}.
\end{proof}

This gives the following corollaries.

\begin{corollary}\label{cor.est}\
We have $\|u_S\|_{H\sp{1-k}(\Omega)} \le C
\|u\|_{H\sp{1-k}(\Omega)}$ for a constant $C$ depending only on
$\Omega$. In particular, $\|u - u_S\|_{H\sp{1-k}(\Omega)} \le C
\|u\|_{H\sp{1-k}(\Omega)}$.
\end{corollary}

\begin{proof}\
The result is well known for $k = 0$ since $u_S$ is the
$B$--orthogonal projection of $u$ onto $S$ (see also C\`ea's
Lemma, \cite{BrennerScott, Ciarlet91}). We shall therefore assume
that $k \ge 1$.

Let $v \in H\sp{k - 1}(\Omega)$ be arbitrary. Let $c \in \CC$ be
such that $\int_{\Omega}(v - c)dx = 0$. Then we can find $V \in
H\sp{k+1}(\Omega)$ such that
\begin{equation*}
    -\Delta V = v - c,\ \int_{\Omega}V dx = 0,\
    \pa_\nu V = 0, \ \text{ and }\
    \|V\|_{H\sp{k + 1}(\Omega)} \le C
    \|v\|_{H\sp{k - 1}(\Omega)},
\end{equation*}
where $C$ is a constant depending only on $\Omega$. Also, chose $w
\in S$ such that
\begin{equation*}
    \|w\|_{H\sp{k + 1}(\Omega)} \le C \|V\|_{H\sp{k + 1}(\Omega)}
    \;\text{ and }\; \|V - w\|_{H\sp{1}(\Omega)}
    \le C h^k \|V\|_{H\sp{k + 1}(\Omega)}.
\end{equation*}
This is possible by Theorem \ref{thm.approx}. Then
\begin{multline*}
    \pairing{u_S}{v} = \pairing{u_S}{v - c}
    = - \pairing{u_S}{\Delta V}
    = - \pairing{u_S}{\Delta V} + \pairing{u_S\vert_{\pa
    \Omega}}{\pa_\nu V} = B(u_S, V) \\
    = B(u_S, w) + B(u_S, V - w)
    = B(u, w) + B(u_S, V - w).
\end{multline*}
Using also Lemmata \ref{lemma.cont.Bu} and \ref{lemma.2}, this
gives
\begin{multline*}
    |\pairing{u_S}{v}| \le C\normH{u}{1-k} \normH{w}{1+k} +
    \normH{u_S}{1}\normH{V - w}{1} \\
    \le C\normH{u}{1-k} \normH{V}{k+1} + C
    h\sp{-k}\normH{u}{1-k} h\sp{k}\normH{V}{k+1}\\
    \le C\normH{u}{1-k} \normH{V}{k+1}
    \le C\normH{u}{1-k} \normH{v}{k-1}.
\end{multline*}
This gives the result since $\normH{u_S}{1-k} := \sup
|\pairing{u_S}{v}|/\normH{v}{k-1}$, $v \neq 0$.
\end{proof}

Similarly,

\begin{corollary}\label{cor.est.b}\
We have $\|u_S\vert_{\pa \Omega}\|_{H\sp{1/2-k}(\pa \Omega)} \le C
\|u\|_{H\sp{1-k}(\Omega)}$ for a constant $C$ depending only on
$\Omega$. In particular, $\|(u - u_S)\vert_{\pa
\Omega}\|_{H\sp{1/2-k}( \pa \Omega)} \le C
\|u\|_{H\sp{1-k}(\Omega)}$.
\end{corollary}

\begin{proof}\ The proof is similar to that of the previous
corollary. Let $v \in H\sp{-1/2 + k}(\pa \Omega)$ be arbitrary.
Let $c \in \CC$ be a constant such that $\int_{\pa \Omega} v  dS =
\int_{\Omega} c dx$. Then we can find a unique $W \in H\sp{1 +
k}(\Omega)$ satisfying
\begin{equation*}
    \Delta W = c,\ \int_{\Omega} W dx = 0,\ \pa_\nu W = v,\
    \text{ and } \normH{W}{k+1} \le C
    \|v\|_{H\sp{-1/2+k}(\pa \Omega)},
\end{equation*}
for a constant $C > 0$ depending only on $\Omega$. Using also
Theorem \ref{thm.approx}, we choose $w \in S$ such that
$\normH{w}{k+1} \le C \normH{W}{k+1}$ and $\normH{W - w}{1} \le
Ch\sp{k} \normH{W}{k+1}$.

Then, using also $\pairing{u_S}{1} = 0$, we obtain
\begin{multline*}
    \pairing{u_S\vert_{\pa \Omega}}{v} =
    \pairing{u_S\vert_{\pa \Omega}}{\pa_\nu W} =
    \pairing{u_S}{\Delta W}  + B(u_S, W) = B(u_S, W) \\
    = B(u_S, w) + B(u_S, W - w) = B(u, w) + B(u_S, W - w).
\end{multline*}
Using Lemmata \ref{lemma.cont.Bu} and \ref{lemma.2}, we then
obtain
\begin{multline*}
    |\pairing{u_S\vert_{\pa \Omega}}{v}| \le C
    \normH{u}{1-k}\normH{w}{1+k} + \normH{u_S}{1}
    \normH{W - w}{1} \\ \le C \normH{u}{1-k}\normH{W}{1+k}
    + C h\sp{-k}\normH{u}{1-k} h\sp{k} \|W\|_{H\sp{1+k}(\Omega)}
    \\ \le C \normH{u}{1-k}  \|W\|_{H\sp{1+k}(\Omega)}
    \le C \normH{u}{1-k}  \|v\|_{H\sp{-1/2+k}(\pa \Omega)},
\end{multline*}
which completes the proof in view of the definition of
$\normH{u_S}{1/2-k}$.
\end{proof}

We complete our sequence of estimates with the following result.

\begin{proposition}\label{prop.f.est}\
Let $k + \gamma \le m + 1$, $k, \gamma \in \ZZ_+$. Assume that
each local approximation space $\Psi_j$ contains the polynomials
of degree $1 + k + \gamma$. Then the error $u - u_S$ satisfies
\begin{equation*}
    \normH{u - u_S}{1 - k - \gamma} \le C h\sp{\gamma}
    \normH{u}{1-k},
\end{equation*}
with a constant $C$ independent of $h$ or the GFEM--space $S$, but
possibly depending on $\Omega$ and the parameters $A$, $B$, $C_j$,
$\kappa$, $\sigma$, and $m$.
\end{proposition}

\begin{proof}\ Let $v \in H\sp{-1 + k + \gamma}(\Omega)$ be
arbitrary. Let $c$ be a constant such that $\pairing{v - c}{1} =
0$. Then there exists a unique $V \in H\sp{1 + k +
\gamma}(\Omega)$ such that
\begin{equation*}
    -\Delta V = v - c,\ \int_\Omega Vdx = 0,\
    \pa_\nu V = 0,\ \text{ and }\  \|V\|_{H^{1+k+\gamma}(\Omega)}
    \le C\|v\|_{H^{-1+k+\gamma}(\Omega)},
\end{equation*}
for a constant $C > 0$ depending only on $\Omega$.

Then, for any $w \in S$,
\begin{multline}\label{eq.EQ}
    \pairing{u - u_S}{v} = \pairing{u - u_S}{v - c} =
    - \pairing{u-u_S}{\Delta V} \\ = - \pairing{u-u_S}{\Delta V}
    + \pairing{(u - u_S)\vert _{\pa \Omega}}{\pa_\nu V}
    = B(u - u_S, V) = B(u - u_S, V - w) \\ =
    -\pairing{u - u_S}{\Delta(V - w)} + \pairing{(u - u_S)\vert_{\pa
    \Omega}}{ \pa_\nu w}.
\end{multline}
Using Theorem \ref{thm.approx}, we chose $w \in S$ such that
$\normH{w}{1 + k + \gamma} \le C \normH{V}{1 + k + \gamma}$ and
$\normH{V - w}{1 + k} \le Ch^{\gamma} \normH{V}{1 + k + \gamma}$.
In particular,
\begin{equation*}
    \|\pa_\nu w\|_{H\sp{-1/2 + k}(\pa \Omega)} = \|\pa_\nu (V -
    w)\|_{H\sp{-1/2 + k}(\pa \Omega)} \le h^{\gamma}\normH{V}{1 + k +
    \gamma}.
\end{equation*}
From $\normH{V}{1+ k + \gamma} \le C
\|v\|_{H\sp{-1+k+\gamma}(\Omega)}$, Corollaries \ref{cor.est} and
\ref{cor.est.b}, and Equation \eqref{eq.EQ}, we then obtain,
\begin{multline*}
    |\pairing{u - u_S}{v}| \le \normH{u - u_S}{1-k}
    \normH{V - w}{1 + k} + \normH{u - u_S}{1-k}
    \|\pa_\nu w\|_{H\sp{-1/2 + k}(\pa \Omega)}\\
    \le C h\sp{\gamma} \normH{u}{1-k}
    \|v\|_{H\sp{-1+k+\gamma}(\Omega)}.
\end{multline*}
The proof is complete.
\end{proof}

Theorem \ref{thm.9.2} and Proposition \ref{prop.f.est} then give
the following result, which is the main result of this paper.

Recall, for the following theorem, that $S = S_{GFEM}$ is the
$GFEM$--space associated to a partition of unity satisfying
Assumptions A--D. Also, recall that we have fixed $u \in
H^{1-k}(\Omega)$ satisfying $\pairing{u}{\Delta \phi} = 0$ for all
$\phi \in \CIc(\Omega)$ and that $u_S$ is the GFEM--approximation
of $u$ (\ie given by Lemma \ref{lemma.disc}).

\begin{proposition}\label{prop.main}\
Assume the local approximation spaces $\Psi_j$ contain the
polynomials of degree $k + \gamma + 1$ and let $A \Subset B
\Subset \Omega$ be admissible open subsets. Then for any $-1 + k
\le -1 + k + \gamma \le m$ and any $l \geq 1$, $k, \gamma \in
\ZZ_+$, we have
\begin{equation*}
    \|u - u_S\|_{H\sp{1}(A)} \le C h\sp{l}
    \|u\|_{H\sp{l + 1}(B)} + C h\sp{\gamma}
    \|u- u_S \|_{H\sp{-1 + k + \gamma}(B)}.
\end{equation*}
\end{proposition}

\begin{proof}\
This follows from Theorem \ref{thm.9.2} and from
\begin{equation*}
    \inf_{\chi \in S} \|u - \chi\|_{H\sp{1}(B)} \le
    Ch\sp{l}\|u\|_{H\sp{l + 1}(B)}.
\end{equation*}
\end{proof}

By taking $l = \gamma$ and using also Proposition \ref{prop.f.est}
and Equation \eqref{eq.int.est}, we obtain

\begin{theorem}\label{thm.main}\
Assume the local approximation spaces $\Psi_j$ contain the
polynomials of degree $k + \gamma + 1$ and let $A_0\Subset \Omega$
be an admissible open subset. Then for any $-1 + k \le -1 + k +
\gamma \le m$, $k, \gamma \in \ZZ_+$, we have
\begin{equation*}
    \|u - u_S\|_{H\sp{1}(A_0)} \le C h^\gamma
    \|u\|_{H\sp{1-k}(\Omega)}.
\end{equation*}
The constant $C$ above is independent of $h$ and $S$, but may
depend on the parameters $A$, $B$, $C_j$, $\kappa$, $\sigma$, $l$,
$\gamma$, and $m$, as well as on $\theta$, the distance between
$\pa A_0$ and $\pa \Omega$.
\end{theorem}

\section{Polynomial local approximation spaces\label{sec.Ver}}

In this section we shall verify that the Assumptions A--D are
verified if we choose $\Psi_j$ to be the space of (restrictions to
$\omega_j$ of) polynomials of degree $\le m$, $m \ge 1$. Most
results of this section are either elementary or well known. We
include them nevertheless for the benefit of the reader and for
completeness.

In this section, the set of polynomials of degree $m$ will be
denoted $\maQ_m$. Also, for any ball $B$ of radius $r$, we shall
denote by $tB$ the ball with the same center as $B$ and radius
$tr$.

\begin{lemma}\label{lemma.T2T}\ There exists a constant $C > 0$,
depending only on $n$, $m$, and $t > 0$, such that for any ball $B
\subset \RR^n$ and any $Q \in \maQ_m$, we have $\|Q\|_{L^2(tB)} \le C
\|Q\|_{L^2(B)}$.
\end{lemma}

\begin{proof}\  For any fixed $B$, $Q \mapsto \|Q\|_{L^2(tB)}$ and $P
\mapsto \|Q\|_{L^2(B)}$ are two norms on the finite dimensional space
$\maQ_m$ of polynomials of degree $\le m$, and hence they are
equivalent. This gives the result, except the independence of $C$ on
$B$. But all balls are affine equivalent and the $L^2$-norm is scaled
by the (square root of the) determinant of the matrix of the affine
transformation.  Thus the constant $C$ can be chosen to be the same
for all balls $B$.
\end{proof}

This gives immediately the following corollary.

\begin{corollary}\label{cor.T2T}\
There exists a constant $C > 0$, depending only on $t > 0$, $n$,
and $m$, such that for any ball $B \subset \RR^n$ and any
polynomial $Q \in \maQ_m$, we have $|Q|_{H^l(tB)} \le C
|Q|_{H^l(B)}$ and $\|Q\|_{H^l(tB)} \le C \|Q\|_{H^l(B)}$, $0 \le l
\le m$.
\end{corollary}

\begin{proof}\ Use Lemma \ref{lemma.T2T} for all derivatives
$Q^{(\alpha)}$, where $\alpha$ is a multi-index such that
$|\alpha| = l$ or $|\alpha| \le l$.
\end{proof}

We now establish to the following ``inverse property.''

\begin{lemma}\label{lemma.inverse}\
There exists a constant $C > 0$, depending only on $m$, $\alpha$, and
$n$ such that $\|Q^{(\alpha)}\|_{L^2(B)} \le C
r^{l-|\alpha|}\|Q\|_{H^l(B)}$ for any $l \le |\alpha| \le m$, any $Q
\in \maQ_m$, and any ball $B$ of radius $r$.
\end{lemma}

\begin{proof}\
Let us prove first the result for $l = 0$.  That is, we need to prove
that $\|Q^{(\alpha)}\|_{L^2(B)} \le C r^{-|\alpha|}
\|Q\|_{L^2(B)}$.

Let $B_1 = B_1(0)$ be the unit ball centered at $0$. Then $Q
\mapsto \|Q^{(\alpha)}\|_{L^2(B_1)}$ is a semi-norm on $\maQ_m$,
the space of polynomials of degree at most $m$, and hence it is
bounded by the norm $Q \mapsto \|Q\|_{L^2(B_1)}$. Thus
$\|Q^{(\alpha)}\|_{L^2(B_1)} \le C_1 \|Q\|_{L^2(B_1)}$. Let $L$ be
an affine transformation mapping $B_1$ onto the ball $B$ of radius
$r$ consisting of the composition of a translation and a dilation
of ratio $r$. Then
\begin{multline*}
    \|Q^{(\alpha)}\|_{L^2(B)} = \det(L)^{1/2} \|Q^{(\alpha)}
    \circ L\|_{L^2(B_1)} = \det(L)^{1/2} r^{-|\alpha|}
    \| (Q \circ L)^{(\alpha)} \|_{L^2(B_1)}\\
    \le C_1 \det(L)^{1/2} r^{-|\alpha|}
    \| Q \circ L \|_{L^2(B_1)} = C r^{-|\alpha|}
    \| Q \|_{L^2(B)},
\end{multline*}
for any $Q \in \maQ_m$.

Assume now that $|\alpha| \ge l > 0$. Choose $\beta \le \alpha$,
$|\beta| = l$. Then
\begin{equation*}
    \|D^{\alpha - \beta} D^\beta Q\|_{L^2(B)}
    \le C r^{-|\alpha - \beta|}\| D^\beta Q \|_{L^2(B)}
    \le C r^{l-|\alpha|} \|Q\|_{H^l(B)}.
\end{equation*}
This completes the proof.
\end{proof}

The relevant ``inverse property'' now follows.

\begin{proposition}\label{prop.inverse}\
There exists a constant $C > 0$, depending only on $\sigma$, $m$,
$\alpha$, and $n$ such that $\|Q^{(\alpha)}\|_{L^2(\Omega)} \le C
r^{l-|\alpha|}\|Q\|_{H^l(B)}$ for any $l \le |\alpha| \le m$, any $Q
\in \maQ_m$, any ball $B$ of radius $r$, and any convex set $\Omega$
contained in $\sigma^{-1}B$.
\end{proposition}

\begin{proof}\ This follows from  Corollary \ref{cor.T2T}
and Lemma \ref{lemma.inverse}.
\end{proof}

We now prove the following elementary lemma.

\begin{lemma}\label{lemma.ref}\
For any $A \Subset \Omega$ and any $k \in \ZZ_+$, we can construct
admissible open sets $A =: B_0 \Subset B_1 \Subset B_2 \Subset
\ldots \Subset B_k \Subset B_{k+1} := \Omega$ such that $\hat C
dist(\pa B_j,\pa B_j) \ge \theta/k$, where $\theta := dist(\pa A,
\pa \Omega)$, provided that $\hat C h < \theta$, where $\hat C$
depending on $n$ only (in particular, $\hat C$ is independent of
$h$).
\end{lemma}

\begin{proof}\ Take $k = 1$, for simplicity. The general result is
proved similarly or by iterating this case. Let $\hat C = 4$. Let
$U$ be the union of all open sets $\omega_j$ at distance at most
$\theta/4$ from $A$.  Let $J$ be the set of indices $j$ such that
$\phi_j \neq 0$ on $U$ and let $G$ be the set where $\sum_{j \in
J} \phi_j = 1$. We then define $B_1$ to be the interior of $G$.
\end{proof}

\subsection{Partition of Unity}
We show in this section that, for a suitable set $\Omega$, we can
choose a family of partitions of unity $\{\phi_j\}$, subordinated
to the covering $\{\omega_j\}$, with $h \to 0$ but with all the
other constants fixed. The proof of the following theorem is not
constructive. A constructive proof, suitable for numerical
implementation, will be included in a forthcoming paper where we
will also discuss the numerical implementation of the GFEM for
boundary value problems with distributional data.

\begin{theorem}\label{thm.exists}\
Let $\Omega$ be a bounded open set with piecewise $C^1$-boundary.
Then there exist constants $A$, $B$, $C_j$, $\kappa$, $m$, and
$\sigma$, such that, for any small enough $h > 0$, we can construct a
partition of unity $\{\phi_j\}$ subordinated to the covering
$\{\omega_j\}$ and satisfying the Assumptions A--D of Section
\eqref{sec.GFEM} for the given value of the parameters $A$, $B$,
$C_j$, $\kappa$, $m$, $\sigma$, and $h$.
\end{theorem}

\begin{proof}\ Let us first construct the covering
$\{\omega_j\}$ and the subsets $\{\omega^*_j\}$. The index $j$ will
belong to a set of points of $\Omega$, the centers of those balls. (So
$\omega_j$ and $\omega^*_j$ will have the same center, namely $j$.)

Let $\epsilon > 0$ be small enough such that, for any $y \in \Omega$
satisfying $dist(y, \pa \Omega)< \epsilon$, there exists a unique $z
\in \pa \Omega$ with $dist(y, \pa \Omega) = dist(y, z)$. Let $\Gamma_r
\subset \Omega$ be the set of points at distance $r$, $r < \epsilon/2$
to $\pa \Omega$. Then $\Gamma_r$ will be a piecewise $C^1$ curve,
bounding a domain diffeomorphic to $\Omega$. Choose on $\Gamma_r$ a
maximal set of points $X_h$ containing the vertices of $\Gamma_r$ and
at distance at least $r/2$ from one another. The maximality of $X_h$
then guarantees that the distance between any two consequtive points
in $X_h$ is at most $r$. Let
\begin{equation*}
    f(r) := \sup_{y \in \pa \Omega}\, dist(y, \Gamma_r).
\end{equation*}
A geometric argument based on the assumption that $\pa \Omega$ is
$C^1$ then shows that $\lim f(r)/r = \mu$, as $r \to 0$, with
$\mu$ finite. Let then $\sigma^{-1} > \max\{4(\mu + 1) , 8 \}$ and
$h := r \sigma^{-1}/2$. This guarantees that the balls with
centers in $X_h$ and diameter $h$ will cover the region between
$\Gamma$ and $\pa \Omega$. Let $Y_h$ be a maximal subset of
$\Omega$ containing $X_h$ and such that any point in $Y_h$ is at
distance at least $r$ to the boundary $\pa \Omega$ and at distance
at least $r/2$ from any other point in $Y_h$.  The choice of $Y_h$
shows that the balls $\omega^*_j$ of diameter $\sigma h = r/2$ and
center at the points $Y_h$ will be disjoint, whereas the balls
$\omega_j$ with the same centers and diameter $h$ will cover the
interior of $\Gamma$.

Let us now construct the partition of unity $\phi_j$.  Let $1 \ge
\psi(t) \ge 0$, $t \ge 0$ be a smooth function such that $\psi(t)
= 1$ if $t \le \sigma$ and $\psi(t) = 0$ if $t \ge 1$. Also, let
$\zeta(t) \ge 0$ be a smooth function such that $\zeta(t) = 0$ if
$t \le 1$ and $1 \ge \zeta(t) > 0$ if $t > 1$ and $\zeta(t) = 1$
if $t \ge 2$. Define then
\begin{equation}\label{eq.prod}
    \eta_j(x) = \psi( dist(x, j)/h ) \prod_{j' \neq j}
    \zeta( dist(x, j')/h ).
\end{equation}
Let $\eta(x) = \sum_j \eta_j(x)$ and $\phi_j(x) =
\eta_j(x)/\eta(x)$. We observe that the number of factors $\neq 1$
in the definition of $\eta_j$ is bounded by a constant independent
of $x$, $j$, and any of the choices above. Therefore $\{\phi_j\}$
is our desired partition of unity.
\end{proof}

Some assumptions on the domain $\Omega$ in the above theorem are
necessary, as shown by the following remark.

\begin{remark}\label{remark.Ver}\
The non-Lipschitz domain
\begin{equation*}
    \Omega_c := \{(x, y), -x^2 \le y \le x^2,
    x^2 + y^2 \le 1, x \ge 0\}
\end{equation*}
will have no covering $\{\omega_j\}$ satisfying the Assumptions
A--D.
\end{remark}

\bibliographystyle{plain}
\bibliography{gfemd}
\end{document}